\documentclass[sigconf]{acmart}
\usepackage{booktabs}
\usepackage{multirow}
\usepackage{makecell}
\usepackage{subcaption}

\AtBeginDocument{%
  \providecommand\BibTeX{{%
    \normalfont B\kern-0.5em{\scshape i\kern-0.25em b}\kern-0.8em\TeX}}}

\copyrightyear{2020} 
\acmYear{2020} 
\setcopyright{acmlicensed}\acmConference[ICAIF '20]{ACM International Conference on AI in Finance}{October 15--16, 2020}{New York, NY, USA}
\acmBooktitle{ACM International Conference on AI in Finance (ICAIF '20), October 15--16, 2020, New York, NY, USA}
\acmPrice{15.00}
\acmDOI{10.1145/3383455.3422522}
\acmISBN{978-1-4503-7584-9/20/10}

\acmSubmissionID{18}

\begin{document}

\title{Understanding Distributional Ambiguity \\ via Non-robust Chance Constraint}

\author{Shumin Ma}
\affiliation {\institution{City University of Hong Kong}}
  \email{shuminma@cityu.edu.hk}

\author{Cheuk Hang Leung}
\affiliation {\institution{City University of Hong Kong}}
\email{chleung87@cityu.edu.hk}

\author{Qi Wu}{\authornote{The corresponding author.}}
\affiliation{\institution{City University of Hong Kong}
}
\email{qiwu55@cityu.edu.hk}

\author{Wei Liu}
\affiliation{\institution{Tencent}}
\email{wl2223@columbia.edu}

\author{Nanbo Peng}
\affiliation{\institution{JD Digits}}
\email{pengnanbo@jd.com}
\renewcommand{\shortauthors}{Ma et al.}

\begin{abstract}
This paper provides a non-robust interpretation of the distributionally robust optimization (DRO) problem by relating the distributional uncertainties to the chance probabilities. Our analysis allows a decision-maker to interpret the size of the ambiguity set, which is often lack of business meaning, through the chance parameters constraining the objective function. We first show that, for general $\phi$-divergences, a DRO problem is asymptotically equivalent to a class of mean-deviation problems. These mean-deviation problems are not subject to uncertain distributions, and the ambiguity radius in the original DRO problem now plays the role of controlling the risk preference of the decision-maker. We then demonstrate that a DRO problem can be cast as a chance-constrained optimization (CCO) problem when a boundedness constraint is added to the decision variables. Without the boundedness constraint, the CCO problem is shown to perform uniformly better than the DRO problem, irrespective of the radius of the ambiguity set, the choice of the divergence measure, or the tail heaviness of the center distribution. Thanks to our high-order expansion result, a notable feature of our analysis is that it applies to divergence measures that accommodate well heavy tail distributions such as the student $t$-distribution and the lognormal distribution, besides the widely-used Kullback-Leibler (KL) divergence, which requires the distribution of the objective function to be exponentially bounded. Using the portfolio selection problem as an example, our comprehensive testings on multivariate heavy-tail datasets, both synthetic and real-world, shows that this business-interpretation approach is indeed useful and insightful.
\end{abstract}

\ccsdesc[500]{Applied computing~Multi-criterion optimization and decision-making}

\keywords{Distributionally robust optimization, Chance constraint, KL divergence, $\phi$-divergence, Heavy-tail distributions, Portfolio selection}


\maketitle

\section{Introduction}\label{sec:Introduction}

Stochastic optimization is widely used in many machine learning algorithms to optimize the expected performance or loss, e.g., the mean squared error for regressions, or the expected discounted return in the context of reinforcement learning \cite{thomas2019concentration}. A sound machine learning model demands reliable estimates of the data-generating distribution. However, uncertainties of the data distribution could arise in many ways: limited observations in the stationary case, time-varying law in the non-stationary case, or the law is subject to policy intervention due to the treatment effect. In robust statistics, formulating a decision-making problem as a DRO problem is a remedy to address the distributional uncertainties in the data \cite{chen2018robust}. 

A typical DRO formulation adds an extra layer of optimization over a set of possible distributions, called the \textit{ambiguity set}, and optimizes the decision variables in the worst-case distribution. There are mainly three ways in the literature to define the ambiguity set. The first is the geometric approach, which allows the parameters of the chosen distribution to vary within certain geometric regions \cite{kim2014deciphering, zhu2009robust, zhu2009worst} such as boxes, ellipsoids, and polyhedrons, etc. The second approach, known as the \textit{moment-based approach}, constructs the ambiguity set by collecting distributions that share the same moment constraints \cite{delage2010distributionally, scarf1957min, chen2011tight, zymler2013distributionally}. The last one, the \textit{statistical distance approach}, uses divergence measures or difference functions between two probability distributions to define the ambiguity set as a ball of distributions \cite{namkoong2017variance, abadeh2015distributionally, esfahani2018data, chen2018robust}. The radius of the ball is referred to as the ambiguity radius. Among the three, the moment-based and the statistical distance approaches address law uncertainties.
In contrast, the geometric approach only addresses the uncertainties in the parameters of a \textit{a prior} fixed distribution, not in its functional form. It does not help if the correct distribution turns out to be lognormal when you think it is instead normal and fine-tune its mean and variance. However, the cost of advancing from parameter uncertainty to law uncertainty is that you lose the interpretability of the ambiguity set because the parameters characterizing it are non-business quantities.  

This paper provides a solution to address this business-interpretation problem. For business applications, a decision-maker would have a hard time relating, e.g., the radius $0.01$ of a KL ball to e.g., product sales, taxi demands, or portfolio returns. The radius $0.01$ is not related to any measures of the business objective. An unavoidable headache for her is how she should decide the size of the ambiguity set. Our idea is straightforward. We want to translate the impact of the ambiguity radius, which lacks business meanings, to the impact of the chance parameters constraining the objective function, which now allows a decision-maker to enter her preferences directly related to the business objective. Take asset allocation as an example, our solution can tell a portfolio manager that setting the ambiguity radius to $0.01$ would be equivalent to asking the optimization not to let the chances of her portfolio return going below $-13\%$ be higher than $2\%$. In this way, the geometry of the ambiguity set, its radius, is connected directly to her granular preference of the objective, the amount of risks she can tolerate.     

This paper makes two primary technical contributions. 
First, our analysis applies to the heavy-tail distributions (e.g., via the Cressie-Read divergence) \cite{glasserman2014robust}, besides the usual light-tail cases using the KL divergence. Heavy-tail distributions, e.g., the lognormal distribution and the student $t$-distribution, are ubiquitous for many business and finance datasets. A DRO problem with an ambiguity set defined by the KL divergence is solvable, however, only when the distribution of the objective function is exponentially bounded \cite{hu2013kullback}, in which case heavy-tail distributions are excluded. Our analysis extends well to the general $\phi$-divergence family, including KL divergence, Burg entropy, $\chi^{2}$-distance, Hellinger distance, Cressie-Read divergence, etc. \cite{namkoong2017variance, hashimoto2018fairness}. The second contribution of this paper is that we establish two connections between a DRO problem and a CCO problem. The first one is that when a bounded constraint is added to decision variables, a DRO problem can be cast as a CCO problem without distributional uncertainties. The second connection is that, without the boundedness constraint, the CCO problem is shown to perform uniformly better than the DRO problem, irrespective of the radius of the ambiguity set, the choice of the divergence measure, or the tail heaviness of the center distribution.

The rest of the paper is organized as follows. In Section \ref{sec:Background}, we provide some background information and the motivation for the proposed optimization problems. Theoretical analysis of the DRO problem and the CCO problem is provided in Section \ref{sec:Theoretical Results}. Section \ref{sec:Solutions} establishes the connection between the DRO problem and the CCO problem under an explicit formulation of the portfolio selection problem. Section \ref{sec:Simulation Experiment} gives numerical experiments, and Section \ref{sec:Conclusion} concludes our findings from both synthetic and empirical data. Due to the page limits, all the proofs are omitted in the main body; however, they can be readily provided once requested.

\section{Problem setup}\label{sec:Background}

\subsection{Notations} 
Let $\mathbf{r}\in\mathbb{R}^{n}$, an $n$-dimensional real-valued random vector, be the vector of asset returns. And suppose the joint probability distribution of $\mathbf{r}$ is $\mathbb{P}$. Let $\mathbb{P}_{0}$ be the nominal probability distribution of $\mathbf{r}$. Let $\mathbf{x}\in\mathbb{R}^{n}$ be the asset allocation strategy, and $\mathbf{e}\in\mathbb{R}^{n}$ be a vector with all entries equal to $1$. Denote the utility function that is concave in $\mathbf{x}$ and associated with $\mathbf{x}$ and $\mathbf{r}$ by $f(\mathbf{x},\mathbf{r})$. We assume that $\mathbf{x}$ lies in a convex set $\mathcal{X}$ and $\mathbb{P}$ belongs to an ambiguity set $\mathcal{U}$. The expectation and variance of a random variable under $\mathbb{P}$ are represented by $\mathbb{E}_{\mathbb{P}}[\cdot]$ and $\mathbb{V}_{\mathbb{P}}[\cdot]$, respectively.

\begin{definition}\label{def:phi-divergence}
(\text{$\phi$-divergence}) Assume that $\phi(t)$ is convex for $t\geq 0$ and that $\phi(1)=0$. Then the $\phi$-divergence $D(\mathbb{Q}||\mathbb{P})$ between distribution $\mathbb{P}$ and distribution $\mathbb{Q}$ is defined as:
	\begin{equation}\label{eqt:phi divergence formula}
	D(\mathbb{Q}||\mathbb{P}):=\int \phi\left(\frac{\mathrm{d}\mathbb{Q}}{\mathrm{d}\mathbb{P}}\right)\mathrm{d}\mathbb{P}=\mathbb{E}_{\mathbb{P}}\left[\phi\left(\frac{\mathrm{d}\mathbb{Q}}{\mathrm{d}\mathbb{P}}\right)\right]:=\mathbb{E}_{\mathbb{P}}\left[\phi\left(L\right)\right].
	\end{equation}
\end{definition}

The quantity $L$ in Eq. \eqref{eqt:phi divergence formula} is called the Radon Nikodym derivative (or likelihood ratio) such that $L\geq 0$ almost surely and $\mathbb{E}_{\mathbb{P}}\left[L\right]=1$. Notice that, for the Radon-Nikodymm derivative $L$ to exist, $\mathbb{Q}$ must be absolutely continuous w.r.t. $\mathbb{P}$. Given the function $\phi$ for a specific $\phi$-divergence, its \textit{conjugate} $\phi^{*}$ is defined as
$\phi^{*}(s):=\sup_{t\geq 0}\{st-\phi(t)\}$.
Table \ref{table:phi divergence example} lists the two divergences used in this paper. But it should be mentioned that, our interpretation of the ambiguity radius actually applies to all the $\phi$-divergences, including Burg entropy, $J$-divergence, $\chi^{2}$-distance, modified $\chi^{2}$-distance, and Hellinger distance. (For more information about the $\phi$-divergence family, see \cite{ben2013robust}).
\begin{table}
  \caption{The two $\phi$-divergences used in this paper. The KL divergence applies to light-tail distributions, while the Cressie-Read divergence is compatible with heavy-tail distributions.}
  \label{table:phi divergence example}
  \centering
  \begin{tabular}{lll}
    \toprule
	 & Kullback-Leibler &  Cressie-Read\\
    \hline
	$\phi(t)$ & $t\log(t)-t+1$ &  $\frac{1-\theta+\theta t-t^{\theta}}{\theta(1-\theta)},\theta\neq 0,1$ \\
	$\phi^{*}(s)$ & $e^{s}-1$& $\frac{(1-s(1-\theta))^{\frac{\theta}{\theta-1}}}{\theta}-\frac{1}{\theta}, s<\frac{1}{1-\theta}$ \\
    \bottomrule
  \end{tabular}
\end{table}

\subsection{Motivation} 
The goal is to maximize the expected utility over a set of admissible allocation strategies $\mathcal{X}$, namely,
\begin{equation}\label{OptOriginal}
\max_{\mathbf{x}\in\mathcal{X}}\mathbb{E}_{\mathbb{P}}[f(\mathbf{x},\mathbf{r})].
\end{equation}

We introduce the ambiguity set $\mathcal{U}$ centered at the nominal distribution $\mathbb{P}_0$ (also called the center distribution in the following context) and controlled by the radius parameter $\rho>0$, that is, $\mathcal{U}:= \{\mathbb{P}:D(\mathbb{P}||\mathbb{P}_{0})\leq\rho\}$. Thus, the distributionally robust counterpart of problem \eqref{OptOriginal} is:
\begin{equation}\label{eqt:OPT1}
\max_{\mathbf{x}\in\mathcal{X}}\min_{\mathbb{P}\in\mathcal{U}}\mathbb{E}_{\mathbb{P}}[f(\mathbf{x},\mathbf{r})].
\end{equation}

For a decision-maker, the ambiguity radius $\rho$ is critical. One cannot set it too large since the optimal utility decreases in $\rho$. However, if it is too small, one loses the robust protection. There is a trade-off in choosing its magnitude in the financial context. In literature, \cite{pardo2005statistical} presents the characteristics of the $\phi$-divergence between the true distribution $\mathbb{P}$ and the nominal distribution $\mathbb{P}_0$, $D(\mathbb{P}||\mathbb{P}_{0})$. Assuming that $\mathbb{P}$ and $\mathbb{P}_{0}$ belong to the same parameterized distribution family with parameter dimension $d$, and that $\phi$ is twice continuously differentiable in a neighborhood of 1 with $\phi^{(2)}(1)>0$, the normalized estimated $\phi$-divergence $\frac{2N}{\phi^{(2)}(1)}D(\mathbb{P}||\mathbb{P}_{0})$ asymptotically (i.e., for the sample size $N \rightarrow \infty$) follows a $\chi^{2}_d$-distribution. This conclusion thus relates the ambiguity radius $\rho$ to a confidence level at which the true distribution $\mathbb{P}$ falls within the ambiguity set. \cite{blanchet2018distributionally} provides one methodology, under the Markowitz's mean-variance portfolio selection framework, to select the ambiguity radius $\rho$ as the smallest radius such that the true asset allocation strategy is included with a given confidence level.

However, in financial practice with real data, the assumption that the true distribution is in the same parameterized family with the center distribution is too strong. A wrong guess of the nominal distribution may lead to a meaningless confidence level interpretation of the ambiguity radius $\rho$. Since the DRO approach is believed to provide robust protection against distributional uncertainty, we are motivated to connect the robust protection to protection provided by traditional risk measures. In particular, the heavy-tail nature of distributions that we are concerned with reminds us of the tail probability protection, the optimization based on which is known as CCO problems. Specifically, we define the CCO problem as:
\begin{equation}\label{eqt:OPT1_prime}
	\max_{\mathbf{x}\in\mathcal{X}} \mathbb{E}_{\mathbb{P}_{0}}[f(\mathbf{x},\mathbf{r})]
	\quad s.t. \quad Pr_{\sim\mathbb{P}_0}(\mathbf{x}^{T}\mathbf{r}\leq-\delta)\leq\epsilon .
\end{equation}

Here, $\delta > 0 $ characterizes a typical investor's loss threshold and $\epsilon >0$ corresponds to the loss probability. The CCO problem in problem \eqref{eqt:OPT1_prime} shares the same objective function as that of problem \eqref{OptOriginal}. The expectation is taken under the nominal distribution $\mathbb{P}_{0}$, not subject to any distributional robustness (the term "non-robust" in the title originates from here). Compared to problem \eqref{OptOriginal}, the new component is the chance constraint with parameters ($\delta$, $\epsilon$) characterizing an investor's tolerance to losses. 

We would build up a performance-based interpretation of the ambiguity radius $\rho$ through the parameters of the chance-constrained problem. To be specific, if under some ambiguity radius $\rho$ and chance constraint parameters ($\delta$, $\epsilon$), problem \eqref{eqt:OPT1} and problem \eqref{eqt:OPT1_prime} achieve the same optimal value, we would say that the robust protection under the ambiguity radius $\rho$ is similar to that of a tail probability protection. In addition, we would also look into how the choice of the allocation strategy set $\mathcal{X}$ and the tail heaviness of the nominal distribution $\mathbb{P}_0$ affect the interpretation of the ambiguity radius $\rho$, given that $\mathcal{X}$ and $\mathbb{P}_0$ are the shared model settings of the two problems \eqref{eqt:OPT1} and \eqref{eqt:OPT1_prime}.

\section{Analysis of DRO and CCO problems}\label{sec:Theoretical Results}

This section is devoted to the theoretical analysis of problems \eqref{eqt:OPT1} and \eqref{eqt:OPT1_prime}. We show that, for general $\phi$-divergences, problem \eqref{eqt:OPT1} can be reformulated as a class of mean-deviation problems with the investor's risk preference parameter controlled by the ambiguity radius $\rho$. Besides, we provide an approximation framework to solve problem \eqref{eqt:OPT1_prime}.

\subsection{Reformulation of the DRO problem \eqref{eqt:OPT1}}\label{sec:Reformulation of OPT1}

Consider the inner optimization problem in problem \eqref{eqt:OPT1}:
\begin{equation}\label{eqt:OPT1_inner layer}
\min_{\mathbb{P}\in\mathcal{U}}\mathbb{E}_{\mathbb{P}}[f(\mathbf{x},\mathbf{r})].
\end{equation}

The Lagrangian dual to problem \eqref{eqt:OPT1_inner layer} is:
\begin{align*}\label{EqLaDual}
& \sup_{\eta_{1}\in\mathbb{R},\eta_{2}\geq 0}\left\{-\frac{1}{\eta_{2}} \sup_{L}\left\{\mathbb{E}_{\mathbb{P}_{0}}[-\eta_{2}(f(\mathbf{x},\mathbf{r})+\eta_{1})L-\phi(L)]\right\}-\eta_{1}-\frac{\rho}{\eta_{2}}\right\} \nonumber \\
 & = \sup_{\eta_{1}\in\mathbb{R},\eta_{2}\geq 0}\left\{-\frac{1}{\eta_{2}}\mathbb{E}_{\mathbb{P}_{0}}[\phi^*(-\eta_{2}(f(\mathbf{x},\mathbf{r})+\eta_{1}))]-\eta_{1}-\frac{\rho}{\eta_{2}}\right\}. 
\end{align*}

The last equality is derived directly from the definition of the conjugate function of $\phi$-divergence. Difficulty in solving the dual problem lies in the term $\mathbb{E}_{\mathbb{P}_{0}}[\phi^*(-\eta_{2}(f(\mathbf{x},\mathbf{r})+\eta_{1}))]$. We hereby follow the idea in \cite{gotoh2018robust} to express optimization \eqref{eqt:OPT1_inner layer} in terms of \textit{Regular Measure of Deviation}, with results summarized in Theorem \ref{thm:Regular Measure of Deviation}.

\begin{theorem}\label{thm:Regular Measure of Deviation}
Let $\phi$ be a closed proper convex function and $\phi^{*}$ be its corresponding conjugate function, respectively. Suppose that under mild conditions, the strong duality holds. Define the regular measure of deviation 
\begin{equation*}
\begin{aligned}
&\mathcal{D}_{\eta_{2},\phi,\mathbb{P}_{0}}(f(\mathbf{x},\mathbf{r})|\mathbb{E}_{\mathbb{P}_{0}}[f(\mathbf{x},\mathbf{r})])\\
:=&\inf_{\eta_{1}}\left\{\eta_{1}+\frac{1}{\eta_{2}}\mathbb{E}_{\mathbb{P}_{0}}\left[\phi^{*}\left(\eta_{2}(\mathbb{E}_{\mathbb{P}_{0}}[f(\mathbf{x},\mathbf{r})]-f(\mathbf{x},\mathbf{r})-\eta_{1})\right)\right]\right\}.
\end{aligned}
\end{equation*} 
Then, optimization \eqref{eqt:OPT1_inner layer} is equivalent to :
$$\mathbb{E}_{\mathbb{P}_{0}}[f(\mathbf{x},\mathbf{r})]-\inf_{\eta_{2}\geq 0}\left\{\frac{\rho}{\eta_{2}}+\mathcal{D}_{\eta_{2},\phi,\mathbb{P}_0}(f(\mathbf{x},\mathbf{r})|\mathbb{E}_{\mathbb{P}_{0}}[f(\mathbf{x},\mathbf{r})])\right\}.$$
\end{theorem}

Furthermore, the quantity $\mathcal{D}_{\eta_{2},\phi,\mathbb{P}_0}(f(\mathbf{x},\mathbf{r})|\mathbb{E}_{\mathbb{P}_{0}}[f(\mathbf{x},\mathbf{r})])$ can be expanded as a series of terms, the coefficients of which can be computed under the nominal distribution $\mathbb{P}_{0}$. By doing so, we can reformulate the DRO problem \eqref{eqt:OPT1} as a single-layer maximization problem. 

\begin{lemma}\label{thm:Regular Measure of Deviation_Taylor expansion new version}
Suppose that $K$ is an even number, $\phi\in\mathcal{C}^{K+1}$ is a convex function which satisfies $\phi(1)=\phi^{(1)}(1)=0$ and $\phi^{(2)}(1)>0$. Assume that $\mathbb{E}_{\mathbb{P}_0}[X^k]< \infty$ for $k\leq K$ and $X$ is defined as $X:=f(\mathbf{x},\mathbf{r})-\mathbb{E}_{\mathbb{P}_{0}}[f(\mathbf{x},\mathbf{r})]$. Then
\begin{equation}\label{eqt:Regular Measure of Deviation_Taylor expansion new version}
\begin{aligned}
&\mathcal{D}_{\eta_{2},\phi,\mathbb{P}_{0}}(f(\mathbf{x},\mathbf{r})|\mathbb{E}_{\mathbb{P}_{0}}[f(\mathbf{x},\mathbf{r})])\\
 &= \sum_{k=1}^{K-1} b_k \mathbb{E}_{\mathbb{P}_{0}}\left[\left(X+\eta_{1}^*\right)^{k+1}\right] \eta_2^k + o(\eta_{2}^{K-1}),
\end{aligned}
\end{equation}
where $b_k = \frac{(-1)^{k+1}z^{(k)}(0)}{(k+1)!}$, and $\eta_{1}^{*}$ is the optimal solution to 
\begin{displaymath}
\min_{\eta_{1}}\sum_{k=1}^{K-1}  b_k\mathbb{E}_{\mathbb{P}_{0}}\left[\left(X+\eta_{1}\right)^{k+1}\right] \eta_{2}^{k}.
\end{displaymath}
Specifically, $z(\cdot)$ is a function satisfying $z(0)=1$, $z^{(1)}(\cdot) = \frac{1}{\phi^{(2)}(z(\cdot))}$ and $z^{(k)}(\cdot)$ can be obtained recursively for $k\geq 2$.
\end{lemma}

Note that the above expansion applies to general utility functions $f(\mathbf{x},\mathbf{r})$ that are concave in $\mathbf{x}$. More importantly, most of the $\phi$-divergences (KL divergence, Cressie-Read divergence, Burg entropy, $J$-divergence, $\chi^{2}$-distance, modified $\chi^{2}$-distance, and Hellinger distance) satisfy the smoothness conditions. Taking KL and Cressie-Read divergence as example, for $K=4$, we can explicitly solve the terms in Eq. \eqref{eqt:Regular Measure of Deviation_Taylor expansion new version}, as are shown in the following corollary.

\begin{corollary}\label{corollary:Regular Measure of Deviation when n is four}
Consider $K=4$. We have the $4^{th}$ order expansion of $\mathcal{D}_{\eta_{2},\phi,\mathbb{P}_{0}}(f(\mathbf{x},\mathbf{r})|\mathbb{E}_{\mathbb{P}_{0}}[f(\mathbf{x},\mathbf{r})])$:
\begin{equation}\label{eqt:Regular Measure of Deviation when n is four}
\begin{aligned}
&\mathcal{D}_{\eta_{2},\phi,\mathbb{P}_{0}}(f(\mathbf{x},\mathbf{r})|\mathbb{E}_{\mathbb{P}_{0}}[f(\mathbf{x},\mathbf{r})]) \\
&=\sum_{k=1}^{3} b_k \mathbb{E}_{\mathbb{P}_{0}}\left[\left(X+\eta_{1}^*\right)^{k+1}\right] \eta_2^k+o(\eta_{2}^{3}),
\end{aligned}
\end{equation}
with $\eta_{1}^{*}$ being the real root to the 3$^{rd}$ order equation 
\begin{equation*}\label{eqt:Regular Measure of Deviation optimal eta1 when n is four}
\begin{aligned}
\sum_{k=1}^{3} & (k+1) b_k \eta_2^k\cdot \eta_1^k + 12b_3\eta_2^3\mathbb{E}_{\mathbb{P}_{0}}[X^{2}]\cdot\eta_1 \\
& + 4b_3\eta_2^3\mathbb{E}_{\mathbb{P}_{0}}[X^{3}]+ 3b_2\eta_{2}^{2}\mathbb{E}_{\mathbb{P}_{0}}\left[X^{2}\right]=0.
\end{aligned}
\end{equation*}

For KL divergence, the coefficients are $b_{1} = 1/2$, $b_{2} = -1/6$, $b_{3} = 1/24$; for Cressie-Read divergence with $\theta>2$, the coefficients are $b_{1}=1/2$, $b_{2}=(\theta-2)/6$, $b_{3}=(\theta-2)(2\theta-3)/24$.
\end{corollary}

\cite{gotoh2018robust} gives a similar expansion of $\mathcal{D}_{\eta_{2},\phi,\mathbb{P}_0}(f(\mathbf{x},\mathbf{r})|\mathbb{E}_{\mathbb{P}_{0}}[f(\mathbf{x},\mathbf{r})])$ in Proposition 3.5. The main difference between our expansion in Eq. \eqref{eqt:Regular Measure of Deviation when n is four} and their expansion lies in the calculation of $\eta_1^*$. In Eq. \eqref{eqt:Regular Measure of Deviation when n is four}, $\eta_1^*$ is directly solved through the polynomial equation, while in \cite{gotoh2018robust}, $\eta_1^*$ is an approximated function of $\eta_2$. 

In the sequel, we take $K = 2$, consider the $2^{nd}$ order expansion of $\mathcal{D}_{\eta_{2},\phi,\mathbb{P}_{0}}(\mathbf{x}^{T}\mathbf{r}|\mathbf{x}^{T}\mathbf{\mu})$ and ignore the higher order terms, which gives

\begin{equation*}
\begin{aligned}
& \min_{\mathbb{P}\in\mathcal{U}}\mathbb{E}_{\mathbb{P}}[f(\mathbf{x},\mathbf{r})] \\
& \approx\mathbb{E}_{\mathbb{P}_{0}}[f(\mathbf{x},\mathbf{r})]-\inf_{\eta_{2}\geq 0}\left\{\frac{\rho}{\eta_{2}}+\frac{\eta_{2}}{2\phi^{(2)}(1)}\mathbb{V}_{\mathbb{P}_{0}}[f(\mathbf{x},\mathbf{r})]\right\}\\
& = \mathbb{E}_{\mathbb{P}_{0}}[f(\mathbf{x},\mathbf{r})]-\sqrt{\frac{2\rho}{\phi^{(2)}(1)}\mathbb{V}_{\mathbb{P}_{0}}[f(\mathbf{x},\mathbf{r})]}.
\end{aligned}
\end{equation*}

The last equality comes as a result of 
\begin{displaymath}
\inf_{\eta_{2}\geq 0}\left\{\frac{\rho}{\eta_{2}}+\frac{\eta_{2}}{2\phi^{(2)}(1)}\mathbb{V}_{\mathbb{P}_{0}}[f(\mathbf{x},\mathbf{r})]\right\} = \sqrt{\frac{2\rho\mathbb{V}_{\mathbb{P}_{0}}[f(\mathbf{x},\mathbf{r})]}{\phi^{(2)}(1)}},
\end{displaymath}
and the minimum is achieved at $\eta_2 = \sqrt{ \frac{ 2\rho\phi^{(2)}(1) }{\mathbb{V}_{\mathbb{P}_{0}}[f(\mathbf{x},\mathbf{r})]}}$. This suggests, when $\rho$ is small, the optimal Lagrangian multiplier $\eta_2$ is also small and the expansion in \eqref{eqt:Regular Measure of Deviation_Taylor expansion new version} is accurate. By taking $\max_{\mathbf{x}\in\mathcal{X}}$ on both sides, we finally achieve the $2^{nd}$ order reformulation of problem \eqref{eqt:OPT1} in Theorem \ref{corollary:OPT1 and OPT1_prime}.

\begin{theorem}\label{corollary:OPT1 and OPT1_prime}
Suppose that $\phi$ is convex, twice continuously differentiable, and that $\phi(1)=\phi^{(1)}(1)=0$ and $\phi^{(2)}(1)>0$. The DRO problem in problem \eqref{eqt:OPT1} is asymptotically equivalent to a mean-deviation problem:
\begin{equation}\label{eqt:OPT1_equivalent version}
\max_{\mathbf{x}\in\mathcal{X}} \left\{\mathbb{E}_{\mathbb{P}_{0}}[f(\mathbf{x},\mathbf{r})]-\sqrt{\frac{2\rho\mathbb{V}_{\mathbb{P}_{0}}[f(\mathbf{x},\mathbf{r})]}{\phi^{(2)}(1)}}\right\}.
\end{equation}
\end{theorem}

Theorem \ref{corollary:OPT1 and OPT1_prime} tells that the ambiguity radius $\rho$ actually controls the investor's risk preference. 

\subsection{Reformulation of the CCO problem \eqref{eqt:OPT1_prime}}\label{sec:Reformulation of OPT1_prime}

Notice that, the chance constraint in problem \eqref{eqt:OPT1_prime} is in the same form as the definition of Value-at-Risk (VaR), a risk measure that focuses on the probability of losses. This motivates us to reorganize the tail chance constraint in problem \eqref{eqt:OPT1_prime} with VaR. The VaR is defined as the minimal level $\gamma$ such that the probability that the portfolio loss $-\mathbf{x}^{T}\mathbf{r}$ exceeds $\gamma$ is below $\epsilon$:
\begin{equation*}
\mathrm{V}_{\epsilon}(\mathbf{x}):=\inf\{\gamma\in\mathbb{R}:Pr_{\sim\mathbb{P}_0}\{-\mathbf{x}^{T}\mathbf{r}\geq \gamma\}\leq\epsilon\}.
\end{equation*}

The equivalent form of the chance constraint in problem \eqref{eqt:OPT1_prime}: $Pr_{\sim\mathbb{P}_0}\{-\mathbf{x}^{T}\mathbf{r}\geq \delta\}\leq\epsilon$ implies that, $\delta$ is included in the set $\{\gamma\in\mathbb{R}:Pr_{\sim\mathbb{P}_0}\{-\mathbf{x}^{T}\mathbf{r}\geq \gamma\}\leq\epsilon\}$. That is to say, the chance constraint can be reorganized with $\mathrm{V}_{\epsilon}(\mathbf{x})$, namely, 
$$Pr_{\sim\mathbb{P}_0}\{-\mathbf{x}^{T}\mathbf{r}\geq \delta\}\leq\epsilon\Leftrightarrow \mathrm{V}_{\epsilon}(\mathbf{x})\leq\delta.$$
Hence, given $\mathbb{E}_{\mathbb{P}_{0}}[\mathbf{x}^{T}\mathbf{r}]=\mathbf{x}^{T}\mathbf{\mu}$, problem \eqref{eqt:OPT1_prime} can be reformulated as
\begin{equation*}
	\max_{x\in\mathcal{X}} 	\mathbf{x}^{T}\mathbf{\mu}
	\quad s.t. \quad \mathrm{V}_{\epsilon}(\mathbf{x})\leq\delta.
\end{equation*}
If $\mathbb{P}_{0}$ is normal, then the VaR can be expressed as
\begin{equation*}
\mathrm{V}_{\epsilon}(\mathbf{x})=\kappa(\epsilon)\sqrt{\mathbf{x}^{T}\Sigma \mathbf{x}}-\mathbf{x}^{T}\mathbf{\mu},
\end{equation*}
where $\kappa(\epsilon)=-\Phi^{-1}(\epsilon)$ and $\Phi^{-1}(\cdot)$ is the inverse of the cumulative distribution function of the standard normal distribution. If $\mathbb{P}_0$ is a member of general elliptical distribution family, \cite{lesniewski2016asymptotics} gives an asymptotic expansion of $\mathrm{V}_{\epsilon}(\mathbf{x})$, which takes the form $\kappa(\epsilon)\sqrt{\mathbf{x}^{T}\Sigma \mathbf{x}}-\mathbf{x}^{T}\mathbf{\mu}$ asymptotically when $\epsilon\rightarrow 0$. For example, if $\mathbb{P}_{0}$ is a student $t$-distribution with degree of freedom parameter $\nu$, then $\kappa(\epsilon)=D\epsilon^{-\frac{1}{\nu}}$, where $D=\left(\frac{c_{n}\pi^{\frac{n-1}{2}}\Gamma(\frac{\nu+1}{2})}{\nu\Gamma(\frac{\nu+n}{2})}\right)^{\frac{1}{\nu}}$, $c_{n}=\frac{\Gamma(\frac{\nu+n}{2})}{\Gamma(\frac{\nu}{2})}\nu^{\frac{\nu}{2}}\pi^{-\frac{n}{2}}$, and $\Gamma(\cdot)$ refers to the gamma function. For distributions other than elliptical distributions, $\sqrt{\frac{1-\epsilon}{\epsilon}}\sqrt{\mathbf{x}^{T}\Sigma \mathbf{x}}-\mathbf{x}^{T}\mathbf{\mu}$ is proved to be a valid approximation of $\mathrm{V}_{\epsilon}(\mathbf{x})$ \cite{ghaoui2003worst,bonami2009exact}. These in all provide the approximation of problem \eqref{eqt:OPT1_prime} reformulated as
\begin{equation}
	\begin{aligned}\label{eqt:OPT1_prime mean deviation form}
		\max_{x\in\mathcal{X}}  \mathbb{E}_{\mathbb{P}_0}[f(\mathbf{x},\mathbf{r})]
		\quad s.t. \quad \kappa(\epsilon)\sqrt{\mathbf{x}^{T}\Sigma \mathbf{x}}-\mathbf{x}^{T}\mathbf{\mu}\leq\delta.
	\end{aligned}
\end{equation}

With the following lemma, we can verify that problem \eqref{eqt:OPT1_prime mean deviation form} is a convex optimization when $\kappa(\epsilon)>0$. For general feasibility set $\mathcal{X}$, problem \eqref{eqt:OPT1_prime mean deviation form} can always be efficiently solved with second-order cone programming (SOCP).

\begin{lemma}\label{thm:Checking Convex function}
Suppose $a>0$. Then the function $a\sqrt{\mathbf{x}^{T}\Sigma \mathbf{x}}-\mathbf{x}^{T}\mathbf{\mu}$ is a convex function of $\mathbf{x}$.
\end{lemma}

\section{Explicit formulations of portfolio selection}\label{sec:Solutions}

In this section, we propose the explicit formulations for portfolio selection problem with $f(\mathbf{x},\mathbf{r}) = \mathbf{x}^T\mathbf{r}$. It only remains to explicitly specify the set $\mathcal{X}$. We begin with the most simple but basic unbounded set $\mathcal{X}:=\left\{\mathbf{x}\in\mathbb{R}^{n}\mid\mathbf{x}^{T}\mathbf{e}=1\right\}$. We would denote the optimal solution and optimal value to optimization \eqref{eqt:OPT1_equivalent version} by $\mathbf{x}^*$ and $v^*$, respectively. The corresponding optimal solution and optimal value to optimization \eqref{eqt:OPT1_prime mean deviation form} are denoted by $\tilde{\mathbf{x}}^*$ and $\tilde{v}^*$, respectively. Throughout the rest of the paper, we would denote $\mathbb{E}_{\mathbb{P}_{0}}[\mathbf{r}]$ and covariance matrix of $\mathbf{r}$ under $\mathbb{P}_{0}$ by $\mathbf{\mu}$ and $\Sigma$, respectively. Then naturally, $\mathbb{E}_{\mathbb{P}_{0}}[\mathbf{x}^{T}\mathbf{r}] = \mathbf{x}^{T}\mathbf{\mu}$, $\mathbb{V}_{\mathbb{P}_{0}}[\mathbf{x}^{T}\mathbf{r}] = \mathbf{x}^{T}\Sigma \mathbf{x}$, $v^* = \mathbf{x}^{*T}\mathbf{\mu}-\sqrt{\frac{2\rho\mathbf{x}^{*T}\Sigma \mathbf{x}^*}{\phi^{(2)}(1)}}$ and $\tilde{v}^* = \tilde{\mathbf{x}}^{*T}\mathbf{\mu}$. 

Recall that in a convex optimization, any local optimum is also a global optimum. This motivates us to study the optimal solution to problem \eqref{eqt:OPT1_equivalent version}, $\mathbf{x}^*$, and the optimal solution to problem \eqref{eqt:OPT1_prime mean deviation form}, $\tilde{\mathbf{x}}^*$, through KKT conditions. The results for $(\mathbf{x}^*, v^*)$ and $(\tilde{\mathbf{x}}^*,\tilde{v}^*)$ are summarized in Theorem \ref{thm:Solution of OPT1_equivalent version} and Theorem \ref{thm:Solution of OPT1 and OPT1_prime}, respectively.

\begin{theorem}\label{thm:Solution of OPT1_equivalent version}
Suppose $\phi^{(2)}(1)>0$. Define $A:=\mathbf{e}^{T}\Sigma^{-1}\mathbf{e}$, $B:=\mu^{T}\Sigma^{-1}\mathbf{e}$, and $C:=\mu^{T}\Sigma^{-1}\mu$. Then problem \eqref{eqt:OPT1_equivalent version} with $f(\mathbf{x},\mathbf{r}) = \mathbf{x}^T\mathbf{r}$ and the feasibility set $\mathcal{X}=\left\{\mathbf{x}\in\mathbb{R}^{n}\mid\mathbf{x}^{T}\mathbf{e}=1\right\}$ has an optimal solution 
when $\rho>\phi^{(2)}(1)(C-{B^{2}}/{A})/2$. And the optimal solution $\mathbf{x}^{*}$ and optimal value $v^*$ are:
\begin{equation*}
\mathbf{x}^{*}=\frac{\Sigma^{-1}(\mu-\lambda^{*}\mathbf{e})}{B-\lambda^{*} A},
\quad v^*=\lambda^*,
\end{equation*}
where $\lambda^{*}=\frac{B}{A}-\frac{\sqrt{B^{2}-A\left(C-2\rho/{\phi^{(2)}(1)}\right)}}{A}$.
\end{theorem}

\begin{theorem}\label{thm:Solution of OPT1 and OPT1_prime}
Suppose $\kappa(\epsilon)>0$. Problem \eqref{eqt:OPT1_prime mean deviation form} with $f(\mathbf{x},\mathbf{r}) = \mathbf{x}^T\mathbf{r}$ and the feasibility set $\mathcal{X}=\left\{\mathbf{x}\in\mathbb{R}^{n}\mid\mathbf{x}^{T}\mathbf{e}=1\right\}$ has an optimal solution when $(\epsilon, \delta)$ satisfies $C-B^{2}/{A}<(\kappa(\epsilon))^{2}<\delta^{2}A+2\delta B+C$ and $B+\delta A>0$. ($A$, $B$, and $C$ defined in Theorem \ref{thm:Solution of OPT1_equivalent version}.) And the optimal solution $\tilde{\mathbf{x}}^{*}$ and optimal value $\tilde{v}^{*}$ are:
$$\tilde{\mathbf{x}}^{*}=\frac{\Sigma^{-1}[(1+\tilde{\lambda})\mu-\tilde{\theta}\mathbf{e}]}{(1+\tilde{\lambda})B-\tilde{\theta}A}, \quad \tilde{v}^{*}=\tilde{\lambda}\delta+\tilde{\theta},$$
where $\tilde{\lambda} =\frac{\sqrt{AC-B^2}}{A \kappa(\epsilon)^{2}-AC+B^{2}}(\frac{\kappa(\epsilon)(B+A\delta)}{\sqrt{A\delta^{2}+2B\delta+C-\kappa(\epsilon)^{2}}} + \sqrt{AC-B^2})$, and $\tilde{\theta}=\frac{(C+\delta B)(\tilde{\lambda}+1)-\tilde{\lambda}\kappa(\epsilon)^{2}}{B+\delta A}$. 

Furthermore, $\tilde{v}^*\geq v^*$, i.e., problem \eqref{eqt:OPT1_equivalent version} always outperforms problem \eqref{eqt:OPT1_prime mean deviation form}.
\end{theorem}

In Theorem \ref{thm:Solution of OPT1 and OPT1_prime}, we first identify the sufficient conditions of $(\epsilon,\delta)$ for the optimization problem \eqref{eqt:OPT1_prime mean deviation form} to be feasible. The comparison between $\tilde{v}^*$ and $v^*$ shows that the CCO reformulation performs uniformly better than the DRO reformulation. Here it should be mentioned that, the outperformance of problem \eqref{eqt:OPT1_equivalent version} over problem \eqref{eqt:OPT1_prime mean deviation form} is not so obvious. At the first glance, it does seem quite straightforward that the objective function in problem \eqref{eqt:OPT1_prime mean deviation form} is always smaller than that in problem \eqref{eqt:OPT1_equivalent version}. While in fact, rather than comparing $\mathbf{x}^{T}\mathbf{\mu}-\sqrt{\frac{2\rho\mathbf{x}^{T}\Sigma \mathbf{x}}{\phi^{(2)}(1)}}$ and $\mathbf{x}^{T}\mathbf{\mu}$ based on the same asset allocation strategy $\mathbf{x}$, we are comparing the two objective functions based on their respective optimal asset allocation strategies, namely, $\mathbf{x}^{*T}\mathbf{\mu}-\sqrt{\frac{2\rho\mathbf{x}^{*T}\Sigma \mathbf{x}^*}{\phi^{(2)}(1)}}$ vs $\tilde{\mathbf{x}}^{*T}\mathbf{\mu}$.

For more complex sets $\mathcal{X}$, we resort to numerical analysis to investigate interpretation of the ambiguity radius $\rho$ through chance constraint parameters.

\section{Experiments}\label{sec:Simulation Experiment}

Sections \ref{SecAccu} and \ref{SecEpsilonDelta} are based on synthetic data to test the reformulation accuracy of the DRO problem \eqref{eqt:OPT1} and to see how the tail heaviness of the nominal distribution $\mathbb{P}_0$ affects the interpretation of $\rho$. Section \ref{sec:Empirical studies} is devoted to a more detailed understanding of the ambiguity radius $\rho$ based on the empirical daily returns of 4 asset classes. And Section \ref{SecHF} uses intraday 5-minute stock returns to test the value of robust protection in the real portfolio selection problem.

\subsection{Reformulation accuracy of problem \eqref{eqt:OPT1}}\label{SecAccu}

In this section, we numerically test the accuracies of the $2^{nd}$ order and the $4^{th}$ order reformulations with respect to the original robust problem \eqref{eqt:OPT1}. The $\phi$-divergence we take is KL divergence, under which problem \eqref{eqt:OPT1} can be exactly solved. And we take the exact optimal value as a benchmark to compare the $2^{nd}$ order and the $4^{th}$ order reformulations.

Table \ref{Table:Comparison of different orders} records the relative errors (in the $3^{rd}$ \& $4^{th}$ columns) w.r.t. the exact optimal value (the $2^{nd}$ column) under KL divergence. It shows that the higher order improvement is particularly notable when data exhibits a heavier tail. In the case of Cressie-Read divergence, which we do not record in the table due to the page limit, we observe a 50 times improvement: when $\rho$ is set to $0.78$, relative error for the $4^{th}$ order reformulation is $1.53\%$, while it is $56.54\%$ for the $2^{nd}$ order reformulation given that the optimal value is $-0.2787$. Here, we assume the ambiguity set under the KL divergence centers at a six-dimensional multivariate exponential distribution with mean=0.2, std=0.2, skewness=2, and kurtosis=6. We set the dimensions to be i.i.d to see a clean impact from the heavy tail. And the center distribution $\mathbb{P}_{0}$ under Cressie-Read divergence is multivariate $t$. We see that the larger the size of the ambiguity set (i.e., larger $\rho$), the better the improvement of the $4^{th}$ order reformulation. In fact, the error reduction is about 10 folds in this example. However, using the $2^{nd}$ order equivalent formulation is good enough to solve problem \eqref{eqt:OPT1} when $\rho$ is small.
\begin{table}[!htbp]
	\caption{Relative errors of the $4^{th}$ order reformulation and $2^{nd}$ order reformulation w.r.t. the optimal value of problem \eqref{eqt:OPT1}. Ambiguity sets are defined by KL divergence centered at a $6$-$d$ exponential distribution.}\label{Table:Comparison of different orders}
	\centering
	\begin{tabular}{rccc} 
		\toprule
		& & \multicolumn{2}{c}{Relative errors} \\
		\cmidrule(lr){3-4}
		 & Optimal value & $4^{th}$ order & $2^{nd}$ order \\ 
		 \hline
		$\mathbf{\rho} =$ {\bf0.01} & 0.1887  & 0.0002$\%$  & 0.1172$\%$ \\ 
		{\bf 0.02} &  0.1841 & 0.0038$\%$ & 0.2397$\%$ \\
		{\bf 0.03} & 0.1807 & 0.0128 $\%$ &  0.3659$\%$ \\
		 {\bf 0.04} & 0.1778 & 0.0274$\%$ & 0.4951$\%$ \\
		{\bf 0.05} &  0.1753 & 0.0479$\%$ &  0.6270$\%$ \\
		{\bf 0.06} &  0.1730 & 0.0748$\%$ &  0.7613$\%$ \\
		{\bf 0.07} &   0.1710 & 0.1082$\%$ & 0.8979$\%$ \\
		{\bf 0.08} &  0.1691 & 0.1483$\%$ &  1.037$\%$ \\
		 {\bf 0.09} & 0.1673 & 0.1951$\%$ & 1.778$\%$\\ 
		\bottomrule
	\end{tabular}
\end{table}

\subsection{Interpretation of $\rho$ under distributions with different tail heavinesses}\label{SecEpsilonDelta}

\begin{figure}[ht]
  \centering
  \includegraphics[width=.4\textwidth]{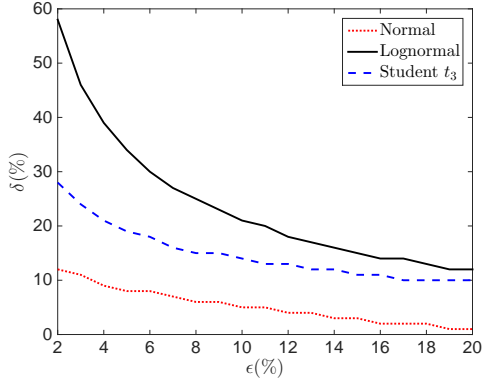}
  \caption{Given $\rho=0.27$, tail heaviness affects the equivalent loss threshold $\delta$.}
  \label{fig:Plot of delta and epsilon for fixed rho}
\end{figure}

This experiment shows that tail heaviness of the nominal distribution $\mathbb{P}_{0}$ indeed affects the interpretation of the ambiguity radius $\rho$. We focus on three distributions for 5 assets: multivariate normal, lognormal distribution and student $t_3-$distribution. The set of allocation strategies is bounded below by -1, and the ambiguity radius $\rho$ is fixed at $0.27$. We plot the results of equivalent $(\epsilon,\delta)$ in Figure \ref{fig:Plot of delta and epsilon for fixed rho}. It shows that, first, the ambiguity radius $\rho$ can be explained by a set of pairs $(\epsilon,\delta)$ in terms of the impact on the optimal value. Second, tail heaviness affects the interpretation of $\rho$ and distributions with heavier tail result in a larger loss threshold for a given loss probability $\epsilon$. 
\vspace{-1em}
\subsection{Empirical studies with daily asset returns}\label{sec:Empirical studies}

To see more clearly the financial interpretation of the ambiguity radius $\rho$, we undergo experiments based on empirical data. We extract past 40 years' daily simple returns of four major asset classes: Equity indexes (DAX,  FTSE, HSI, NASDAQ, NIKKEI250, SP500) , US Treasuries (2year, 10year, 30year), Currencies (AUD, CHF, EUR, GBP, JPY) and Commodities (Crude oil, Silver, Gold). For the DRO problem, we use the Cressie-Read divergence instead of KL divergence since all data exhibits quite heavy tail. For the CCO problem, we choose the negative daily return threshold $-\delta$ to be the $3\%$ empirical quantile of the daily simply return series for each asset class so that they can differ across assets. We choose the chance level $\epsilon$ to be $2\%$ and $5\%$, mincing (rounded) event frequencies at quarterly (4 out of 252) and monthly (12 out of 252) so that investors can relate $\epsilon$ to the degree of event rareness. The portfolio weights are constrained to be bounded below by -1. Both multivariate $t$- and normal distributions are tested as the center $\mathbb{P}_{0}$ of the ambiguity set $\mathcal{U}$ when fitting data. Also, we test both the $4^{th}$ order and $2^{nd}$ order reformulations of the DRO problem.

\begin{table}[!tbp]
    \caption{The equivalent ambiguity radius $\rho$ of the DRO problem for the four asset classes: (\subref{tab:equity}) Equity, (\subref{tab:bond}) US Treasury, (\subref{tab:fx}) Currency, and (\subref{tab:commodity}) Commodity. The loss threshold $\delta$ is taken as the negative value of the $3\%$ empirical quantile of the daily simple return series for each asset class, thus is different across assets. We compare the portfolio performance within each asset class based on the choice of the center distribution (either multivariate student $t$- or normal distributions) under both the $4^{th}$ order and $2^{nd}$ order reformulations of the DRO problem. The percentage number in the round brackets under the equivalent ambiguity radius $\rho$ records the corresponding optimal portfolio annualized return. Bold numbers emphasize the better portfolio return performance at a given pair of $(\epsilon, \delta)$ under a given solution framework of the DRO problem.} 
    \label{Table:Empirical study}
    \centering 
    \begin{subtable}[t]{.45\textwidth}
        \centering
        \caption{Equity: $\delta$ = 3.35$\%$.}
        \label{tab:equity}
        \begin{tabular}{lcccc}
            \toprule
            & \multicolumn{2}{c}{$4^{th}$ order} & \multicolumn{2}{c}{$2^{nd}$ order}\\
            \cmidrule(lr){2-3} \cmidrule(lr){4-5}
            & Student $t$ & Normal & Student $t$ & Normal \\
            \hline
            $\epsilon$ = 2$\%$ & \textbf{3.5e-4} & 1.2e-4  & \textbf{6.1e-4}& 1.2e-4 \\
            & (\textbf{30.7$\%$}) & (15.3$\%$) &  (\textbf{30.7$\%$}) &  (15.3$\%$) \\
            \hline
            $\epsilon$ = 5$\%$ & \textbf{3.4e-4} & 1.2e-4 & \textbf{6.1e-4}  & 1.2e-4 \\
            &  (\textbf{39.2$\%$}) &  (19.8$\%$) & (\textbf{39.2$\%$}) & (19.8$\%$) \\
            \bottomrule
        \end{tabular}
    \end{subtable}
    \\
    ~
    \begin{subtable}[t]{.45\textwidth}
        \centering
        \caption{US Treasury: $\delta$ = 6.58$\%$.}
        \label{tab:bond}
        \begin{tabular}{lcccccccc}
            \toprule
            & \multicolumn{2}{c}{$4^{th}$ order} & \multicolumn{2}{c}{$2^{nd}$ order}\\
            \cmidrule(lr){2-3} \cmidrule(lr){4-5}
            & Student $t$ & Normal & Student $t$ & Normal \\
            \hline
            $\epsilon$ = 2$\%$ & \textbf{2e-6} & 2.8e-14  & \textbf{9.5e-6}& 2.8e-14 \\
            & (\textbf{-1.1$\%$}) & (-2.6$\%$) &  (\textbf{-1.1$\%$}) &  (-2.6$\%$) \\
            \hline
            $\epsilon$ = 5$\%$ & \textbf{2e-6} & 2.8e-14 & \textbf{4.8e-6}  & 2.8e-14 \\
            &  (\textbf{$0.7\%$}) &  (-2.6$\%$) & (\textbf{0.7$\%$}) & (-2.6$\%$) \\
            \bottomrule
        \end{tabular}
    \end{subtable}
    \\
    ~
    \begin{subtable}[t]{.45\textwidth}
        \centering
        \caption{Currency: $\delta$ = 1.40$\%$.}
        \label{tab:fx}
        \begin{tabular}{lcccccccc}
            \toprule
            & \multicolumn{2}{c}{$4^{th}$ order} & \multicolumn{2}{c}{$2^{nd}$ order}\\
            \cmidrule(lr){2-3} \cmidrule(lr){4-5}
            & Student $t$ & Normal & Student $t$ & Normal \\
            \hline
            $\epsilon$ = 2$\%$ & 2.6e-4 & \textbf{6.1e-5}  & 3.1e-4& \textbf{6.1e-5} \\
            & (2.3$\%$) & (\textbf{3.6$\%$}) &  (2.3$\%$) &  (\textbf{3.6$\%$}) \\
            \hline
            $\epsilon$ = 5$\%$ & 1.5e-4 & \textbf{3.1e-5} & 3.1e-4  & \textbf{3.1e-5} \\
            &  (4.4$\%$) &  (\textbf{5.0$\%$}) & (4.4$\%$) & (\textbf{5.0$\%$}) \\
            \bottomrule        
            \end{tabular}
    \end{subtable}
    \\
    ~
    \begin{subtable}[t]{.45\textwidth}
        \centering
        \caption{Commodity: $\delta$ = 4.4$\%$.}
        \label{tab:commodity}
        \begin{tabular}{lcccccccc}
            \toprule
            & \multicolumn{2}{c}{$4^{th}$ order} & \multicolumn{2}{c}{$2^{nd}$ order}\\
            \cmidrule(lr){2-3} \cmidrule(lr){4-5}
            & Student $t$ & Normal & Student $t$ & Normal \\
            \hline
            $\epsilon$ = 2$\%$ & \textbf{9.6e-5} & 3.7e-9  & \textbf{1.5e-4}& 3.7e-9 \\
            & (\textbf{17.3$\%$}) & (4.6$\%$) &  (\textbf{17.3$\%$}) &  (4.6$\%$) \\
            \hline
            $\epsilon$ = 5$\%$ & \textbf{6.5e-5} & 1.9e-9 & \textbf{7.6e-4}  & 1.9e-9 \\
            &  (\textbf{22.7$\%$}) &  (4.6$\%$) & (\textbf{22.6$\%$}) & (4.6$\%$) \\
            \bottomrule
         \end{tabular}
    \end{subtable}
\end{table}

Table \ref{Table:Empirical study} (\subref{tab:equity})-(\subref{tab:commodity}) report the equivalent ambiguity radius $\rho$ of the DRO problem, together with the corresponding optimal portfolio return (annualized), at a given pair of CCO parameters ($\epsilon$, $\delta$) for the four asset classes, respectively. Take Table \ref{Table:Empirical study}(\subref{tab:equity}) as an example. There are 2 rows, 4 columns and 8 entries in total. Each row corresponds to the choice of the parameter $\epsilon$, and each column corresponds to the choice of the reformulation framework of the DRO problem and the choice of the center distribution $\mathbb{P}_0$. The upper number in one entry records the equivalent ambiguity radius $\rho$, while the lower number in the round brackets records the corresponding optimal portfolio annualized return. With other parameter fixed, we compare the optimal portfolio returns between the multivariate student $t$-distribution and the normal distribution, and label the entry numbers with a larger portfolio return in bold black.

We read from Table \ref{Table:Empirical study} that, by relating the size parameter $\rho$ of the ambiguity set in the DRO problem to the CCO chance parameters, it then becomes tangible, without which even the appropriate order is hard to guess. In our tests, its magnitude can range from $10^{-4}$ to $10^{-14}$ depending on asset classes and on the investor's tolerance level. What's more, the heavy-tail nature of financial data demands the usage of divergence measures (e.g., the Cressie-Read divergence) that allow heavy-tail distribution if one takes the robust approach for portfolio optimization. Ambiguity sets constructed by the KL divergence, however, require the objective function to be exponentially bounded, which exclude important heavy-tail distributions used ubiquitously for financial asset returns, e.g., the student $t$-distribution. Among the 16 tests in Table \ref{Table:Empirical study}, the larger return in bold shows 12 favor fitting data with $\mathbb{P}_0$ as multivariate $t$-distributed. 
\vspace{-0.9em}
\subsection{High frequency empirical setting}\label{SecHF}

We collect the intraday 5-minute asset returns of 15 stocks \footnote{The ticker codes for the selected 15 stocks are: 00001, 00005, 00016, 00027, 00388, 00688, 00700, 00883, 00939, 00941, 01299, 01398, 01928, 02318, 03988.} that are selected from the 50 Hang Seng Index constituent stocks based on the market cap and daily turnover. The data spans from Dec 1st, 2014 to Dec 1st, 2017, and consists of roughly 39,390 observations with information of the first and the last half hours in each trading day excluded.  

\begin{figure}[ht]
  \centering
  \includegraphics[width=.4\textwidth]{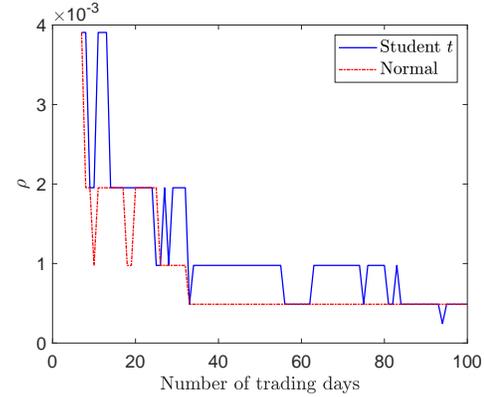}
  \caption{With $(\epsilon,\delta) = (3\%, 0.28\%)$ fixed, the equivalent ambiguity radius $\rho$ goes down and converges as data of more trading days is available.}
  \label{Fig_rho_time}
\end{figure}

The first experiment illustrates the trend of equivalent ambiguity radius $\rho$ as more empirical data is available. As in the last experiment, we use the Cressie-Read divergence and set the loss probability $\epsilon = 3\%$ and $\delta = 0.28\%$ (the $3\%$ empirical quantile of the return series over 100 trading days). The asset allocation strategy is bounded from below by $-1$. We apply the $4^{th}$ order reformulation to solve the DRO problem and test both multivariate $t$- and normal distribution as the nominal distribution $\mathbb{P}_0$. To begin with, we compute the equivalent ambiguity radius $\rho$ based on the first 6 consecutive trading days of 5-minute return series. Then we move forward to include one more trading day's sample data and obtain the next equivalent $\rho$. Figure \ref{Fig_rho_time} plots the series of equivalent $\rho$ with each $\rho$ stamped with how many trading days' data the computation is based on. 

Figure \ref{Fig_rho_time} shows that, to achieve the same level of tail probability protection, the equivalent ambiguity radius $\rho$ goes down and converges as more data is available. Such a conclusion is within expectation because the more available data, the more information and thus fewer uncertainties are over the underlying distribution. The second observation accords with the conclusion in Figure \ref{fig:Plot of delta and epsilon for fixed rho}, that is, even with the same empirical data set, the tail heaviness assumption of the center distribution affects the interpretation of the ambiguity radius $\rho$. Robust portfolio optimization centered with heavy-tail distributions requires a larger range of robust protections to achieve the same tail probability level.

\begin{figure}[ht]
  \centering
  \includegraphics[width=.4\textwidth]{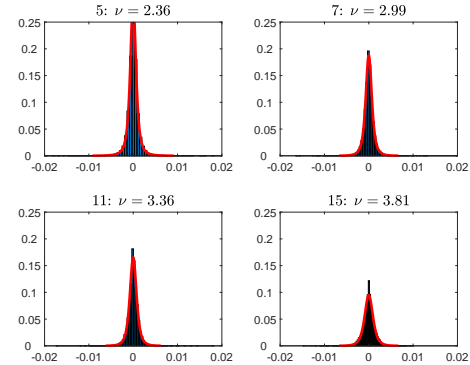}
  \caption{Fitting performance in a student $t$-distribution for stock 5, 7, 11 and 15. The degree of freedom parameter $\nu$ is shown to range from 2.36 to 3.81, which verifies that the distributions of intraday 5-minute returns are indeed heavy-tailed.}
  \label{Fig_tfit}
\end{figure}

Then, we fit the returns of each single stock to a univariate student $t$-distribution to verify that the distribution of high frequency financial data indeed exhibits heavy tail. The degree of freedom parameter, which quantifies the tail heaviness, is shown to range from 2.36 to 3.81 among the 15 stocks. Figure \ref{Fig_tfit} shows the fitting results of 4 stocks accompanied with the degree of freedom parameter $\nu$ in the title position. As it suggests, assuming the nominal distribution of the returns as a student $t$-distribution is rather reasonable.

The last experiment focuses on the value of robust protection in portfolio optimization. In real practice, portfolio optimization under a distributional robust framework is needed to protect investors from uncertainties arising from both the limited historical data and future distributional changes. It is necessary for a trader to frequently rebalance the portfolio to accommodate fluctuations in distributions. As we would demonstrate, the robust protection actually helps improve the portfolio performance, especially when compared with portfolios that are selected either based on the nominal distribution (namely, problem \eqref{OptOriginal}) or under the classical Mean Variance framework. The Mean Variance model we take is:
\begin{equation*}
	\min_{x\in\mathcal{X}} \mathbf{x}^{T}\Sigma \mathbf{x} \quad s.t. \quad \mathbf{x}^{T}\mathbf{\mu} \geq r_{target} .
\end{equation*}
We divide the whole 3-year datasets into two consecutive parts. With the first 2-year data, we fit it to a 15-d student $t$-distribution and establish the equivalent ambiguity radius $\rho$ = 2.4e-4 and optimal return 0.35e-4, given chance constraint parameters $(\epsilon,\delta)$ = $(3\%,39\text{e-4})$. Then with the last-year data as a test set, we backtest the portfolio performance with three asset allocation strategies solved respectively by the DRO problem, the nominal optimization problem, and the Mean Variance problem. For the DRO problem, we fix $\rho$ = 2.4e-4, and for the Mean Variance problem, we fix $r_{target}$ = 0.35e-4. Under each optimization framework, the asset allocation strategy is not constant throughout the whole testing period. We rebalance the portfolio in the frequency of every 5 minutes/hour/half day/day. For each rebalancing, we always use its past 4 months of trading data to solve the optimal allocation strategy and then apply the strategy to next incoming 5 minutes/hour/half day/day. Table \ref{Tb_rebalance} summarizes the statistics of the return series based on different strategies and rebalancing frequencies. 
\begin{table}[!htp]
\caption{Statistics of the 3 return series constructed by 5-minute/hourly/half-day/daily rebalanced allocation strategies solved by the DRO problem, the nominal problem, and the Mean Variance problem, respectively.}
\centering
\begin{tabular}{l ccc}
\toprule
 &   DRO &  Nominal &  Mean Variance \\
 \midrule
\multicolumn{3}{l}{\textbf{5-minute rebalancing}} & \\
Mean (e-4)  &  \textbf{3.68} & 3.24 & 0.77  \\
Variance(e-6)     & 56.5 & 199 & 3.48  \\
Skewness & \textbf{0.66} & 0.17 & 0.61 \\
\midrule
\multicolumn{3}{l}{\textbf{Hourly rebalancing}} & \\
Mean (e-4)  &  \textbf{3.19} & 2.67 & 0.59  \\
Variance(e-6)     & 56.3 & 200 & 3.52  \\
Skewness & \textbf{0.64} & 0.13 & 0.53 \\
\midrule
\multicolumn{3}{l}{\textbf{Half-day rebalancing}} & \\
Mean (e-4)  &  \textbf{2.7} & 2.3 & 0.48  \\
Variance(e-6)     & 55.9 & 198 & 3.47  \\
Skewness & \textbf{0.62} & 0.09 & 0.43 \\
\midrule
\multicolumn{3}{l}{\textbf{Daily rebalancing}} & \\
Mean (e-4)   & \textbf{1.69} & 1.0 & 0.24  \\
Variance(e-6) & 55.6 & 190 & 3.43 \\
Skewness  & \textbf{0.66} & 0.054 & 0.52 \\
\bottomrule
\end{tabular}
\label{Tb_rebalance}
\end{table}%

Table \ref{Tb_rebalance} shows that, the dynamic allocation strategy under a robust framework always outperforms that without a robust protection and the classical Mean Variance strategy. The outperformance can be at most 7 times, depending on the rebalancing frequency. And the DRO strategy keeps a medium level of volatility, neither too aggressive nor too conservative to gain low returns. What's more, the highest skewness for the DRO strategy also highlights its inclination to more gains than losses. Last but not the least, although the outperformance of a DRO strategy is consistent between different rebalancing frequencies, an investor benefits from more frequent rebalancing with returns far more than doubled under whatever portfolio selection framework. 
\vspace{-0.9em}
\section{Conclusions}\label{sec:Conclusion}

We delved into the ambiguity radius for DRO problems with a distributional ambiguity set defined by $\phi$-divergence. We showed that for general $\phi$-divergences, a DRO optimization problem is asymptotically equivalent to a mean-deviation problem, where the risk preference parameter is controlled by the ambiguity radius. We used a portfolio selection example to demonstrate that, when the investment strategy is bounded, the ambiguity radius can be cast as a chance constraint in a deterministic optimization with the same objective. Otherwise, within the set of unbounded investment strategies, a chance-constrained deterministic optimization consistently performs better than the DRO problem. Through extensive experiments with both synthetic and empirical data, we concluded that, to achieve the same level of tail probability protection, a DRO problem centered at heavy-tail distributions requires a larger ambiguity set.

\section*{Acknowledgments}
Qi WU acknowledges the GRF support from the Hong Kong Research Grants Council under 14211316 and 14206117.

%
\bibliographystyle{ACM-Reference-Format}
\bibliography{ICAIF2020_conference}

\appendix

\end{document}